\documentclass[a4paper,12pt]{article}
\usepackage{amsfonts,amssymb,pdfpages, xcolor}

\title{ Nonlinear Generalized Functions: their origin, some developments and recent advances}
\author{
J.F. Colombeau*, \\Institut Fourier, Universit\'e de Grenoble (retired).}
\date{}
\begin{document}
\maketitle

\begin {abstract}  We expose some simple facts at the interplay between mathematics and the real world, putting in evidence mathematical objects " nonlinear generalized functions" that are needed to model the real world, which appear to have been generally  neglected up to now by  mathematicians. Then we describe how a "nonlinear theory of generalized functions" was obtained inside the Leopoldo Nachbin group of infinite dimensional holomorphy between 1980 and 1985 **. This new theory permits to multiply arbitrary distributions and contains the above mathematical objects,  which shows that the features of this theory  are natural and unavoidable for a mathematical description of the real world. Finally we present direct applications of the theory such as existence-uniqueness for systems of PDEs without classical solutions and calculations of shock waves for systems in non-divergence form  done between 1985 and 1995  ***, for which we give    three examples of different nature (elasticity, cosmology, multifluid flows).

\end {abstract} 

\textit{*  work  done under  support of FAPESP, processo 2011/12532-1, and thanks to the hospitality of the IME-USP.}\\
\textit{ ** various supports from Fapesp, Finep and Cnpq between 1978 and 1984.}\\
\textit{***   support of "Direction des Recherches, \'Etudes et Techniques", Minist\`ere de la D\'efense, France,   between 1985 and 1995.}\\
\textit{jf.colombeau@wanadoo.fr}\\

\textbf{1. Mathematics and the real world.} Let $H$ be the Heaviside function. It is a function equal to 0 if $x<0$, to 1 if $x>0$ and not defined if $x=0$. For $p,q\in \mathbb{N}$ let us compute the integral
\begin{equation} I=\int_{-\infty}^{+\infty}(H^p(x)-H^q(x))H'(x)dx. \end{equation}
Physicists solve the problem with formal calculations:
$$ I=[\frac{H^{p+1}}{p+1}]_{-\infty}^{+\infty}-[\frac{H^{q+1}}{q+1}]_{-\infty}^{+\infty}=\frac{1}{p+1}-\frac{1}{q+1} \not=0.$$
The justification is that the Heaviside function $H$ is an idealization of a smooth function whose values  jump from 0 to 1 on a very small interval: before idealization the calculations on the smooth function are well defined
and the result should of course be retained after idealization. Idealization is needed for simplification in the manipulation of the mathematical models of physics. Therefore the concept of Heaviside function is needed but one is forced to consider that $H^p \not=H^q$ if $p\not= q$.\\

In Lebesgue integration theory, therefore in Schwartz distribution theory, one has $ H^p=H^q$ even if $p\not= q$ since they differ only at the point $x=0$ which has Lebesgue measure 0. But one cannot multiply by $H'$ and integrate as requested in (1).\\

To be more convinced is it possible to visualize different Heaviside functions from direct observation of the real world? The answer is yes: it suffices to have in mind compound shock waves which you can produce with a solid metallic ruler and a hammer if you knock strongly enough (figure 1). \\

The shock wave that propagates in the ruler is made of an elastic region followed by a plastic region. In the elastic region the "stress" $S$ of the electrons increases but they remain able to maintain the cristalline structure of the metal. The plastic region begins when the electrons become unable to maintain this structure: then in absence of their action the solid behaves like a fluid with definitive  deformation and possible breakings.  Therefore if the physical variables density and velocity vary throughout the shock wave (along both regions) the stress varies only in the elastic region (in which it increases in absolute value) and then remains fixed at its maximum (in absolute value) in the plastic region.
A 1-D model is provided by the set of equations
\begin{equation}\rho_t+(\rho u)_x=0,\end{equation}
\begin{equation}(\rho u)_t+(\rho u^2)_x+(p-S)_x=0,\end{equation}
\begin{equation}S_t+uS_x-k^2(S)u_x=0, k^2(S)=1 \ if |S|<S_0, k^2(S)=0 \  if \ |S|\geq S_0,
\end{equation}
\begin{equation}p=const.\rho,\end{equation}


where $\rho$=density, $u$=velocity, $p$=pressure, $S$=stress deviation tensor. $|S|$ increases in the elastic region; it stops when it reaches the value $S_0$: then the material  is in the plastic region. (2) and (3) are respectively the continuity and the Euler equations; (4) is a differential form of Hooke's law in which passage from the elastic region to the plastic region is modeled in the Lam\'e coefficient $k^2$, (5) is some simplified state law linking pressure and density.\\ 
\\
\\
\\
\\
\\
\\
\\
\\
\\
\\
\\
\\
\\
\\
\\
\\
\\
\\
\\
\\
\\
\\
\\
\\
\\
\\
\\
\\
\\
\\
\\
\\
\\
\\
\\
\\

 \begin{figure}[h]
\centering \includepdf[width=\textwidth]{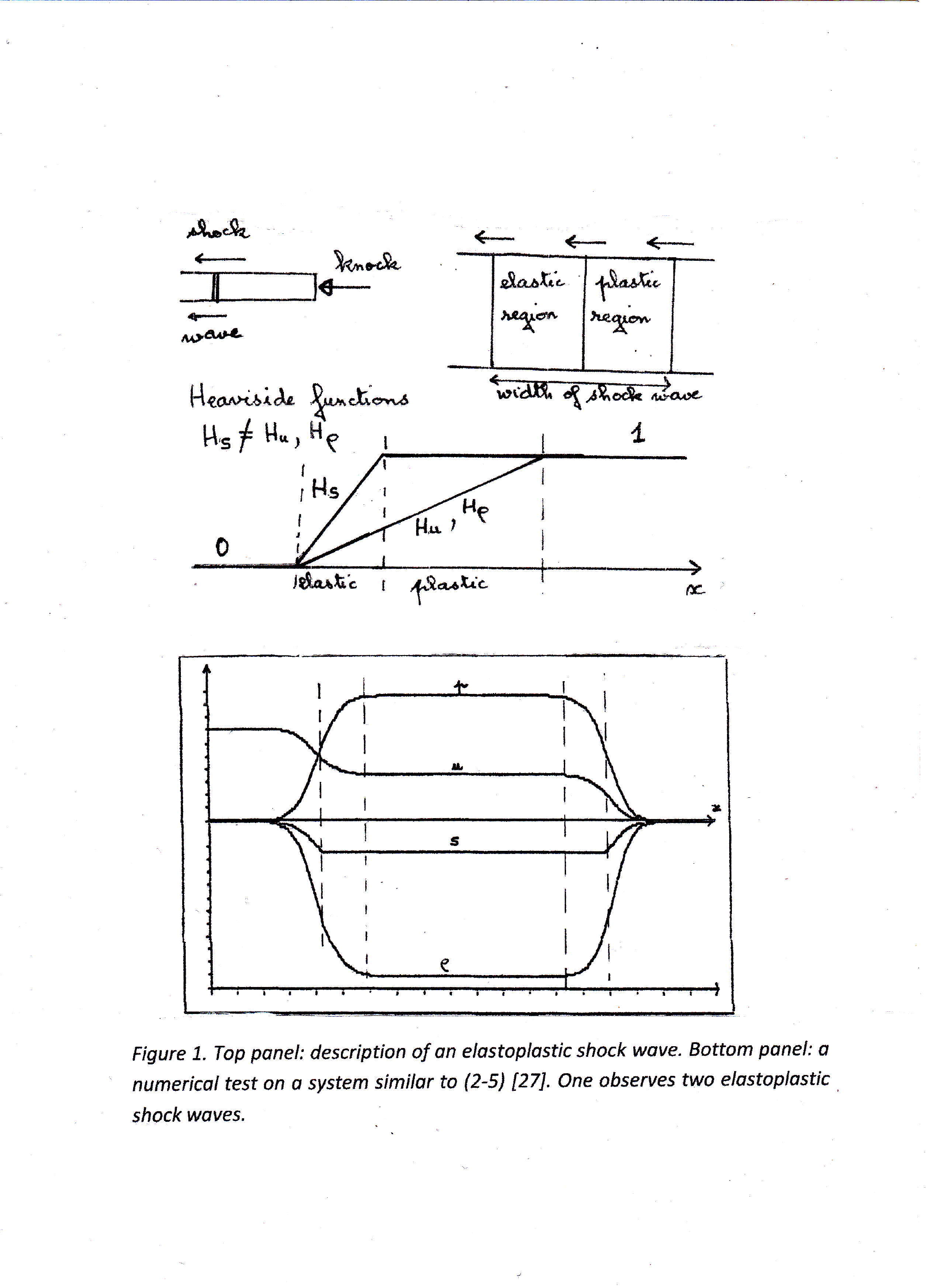}
\end{figure}
The integral $\int u(x,t) S_x(x,t) dx$ (used to try to compute a jump condition for equation 4) has very different values depending on the Heaviside functions representing $u$ and $S$. For very strong shocks the elastic region is  small relatively to the plastic region and therefore the Dirac delta distribution $S_x$ is located on the " foot" of the Heaviside function that represents the jump of $u$, therefore the  integral $\int u(x,t) S_x(x,t) dx$ is close to 0 if $u$=0 on this foot. On the other hand for weak shocks the plastic region can be void: the shock is purely elastic and the variables $\rho, u, S$ vary throughout the shock wave: admitting roughly they  are represented by the same Heaviside function $H$ the integral   $\int u(x,t) S_x(x,t) dx$ has a  far from zero value. See \cite {Colombeau4, Bialivre} for details.\\

Here is an example of another idealization that does not make sense in distribution theory. It appears in form of a square root of the Dirac delta distribution.  The Dirac delta distribution can be asymptotically represented by a function equal to $\frac{1}{\epsilon}$ on $[-\frac{\epsilon}{2},\frac{\epsilon}{2}]$ , 0 elsewhere; what about the object asymptotically represented by the square root of this object? It should be denoted by $\sqrt\delta$. In the viewpoint  of Schwartz distribution theory such an object should be 0 because the integral limited by its graph and the $x$-axis is equal to $\sqrt\epsilon \rightarrow 0$. Does such an object  $\sqrt\delta$ exist in the real world?\\

In general Relativity it was discovered from formal calculations that the Reissner-Nordstrom field of a charged ultra-relastivistic (very fast rotating) black hole is null while the  energy of this (null) field is not null!  R. Steinbauer \cite{Steinbauer} solved this paradox by discovering that the field itself is a mathematical object such as $\sqrt\delta$ above, which looks null in distribution theory but whose square
 (the energy of the field) is similar to a Dirac delta distribution i.e. not at all null.\\

The Keyfitz-Kranzer system of singular shocks
\begin{equation}u_t+(u^2-v)_x=0,\end{equation}
\begin{equation}v_t+(\frac{u^3}{3}-u)_x=0,\end{equation}

has been studied in \cite{Keyfitz1,Keyfitz2, Sanders} and  the observed solutions from numerical tests are given in \cite{Keyfitz2, Sanders}. Here is a depiction of a singular shock wave:\\

In \cite{Keyfitz1} p.426, explicit solutions are found in the form: $\xi=\frac{x}{t}, s=$ shock speed,
$$u(x,t,\epsilon)= u_0+(u_1-u_0)h(\frac{\xi-s}{\epsilon^p}) +\frac{a}{\sqrt{\epsilon}} \rho(\frac{\xi-s}{\epsilon}),$$
$$v(x,t,\epsilon)= v_0+(v_1-v_0)h(\frac{\xi-s}{\epsilon^p}) +\frac{a^2}{\epsilon} \rho^2(\frac{\xi-s}{\epsilon}),$$
\\
\\
\\
\\
\\
\\
\\
\\
\\
\\
\\
\\
\\
\\
\\
\\
\\
\\
\\
\\
\\
\\
\\
\\
\\
\\
\\
\\
\\
\\
\\
\\
\\
\begin{figure}[h]
\centering
\includepdf[width=\textwidth]{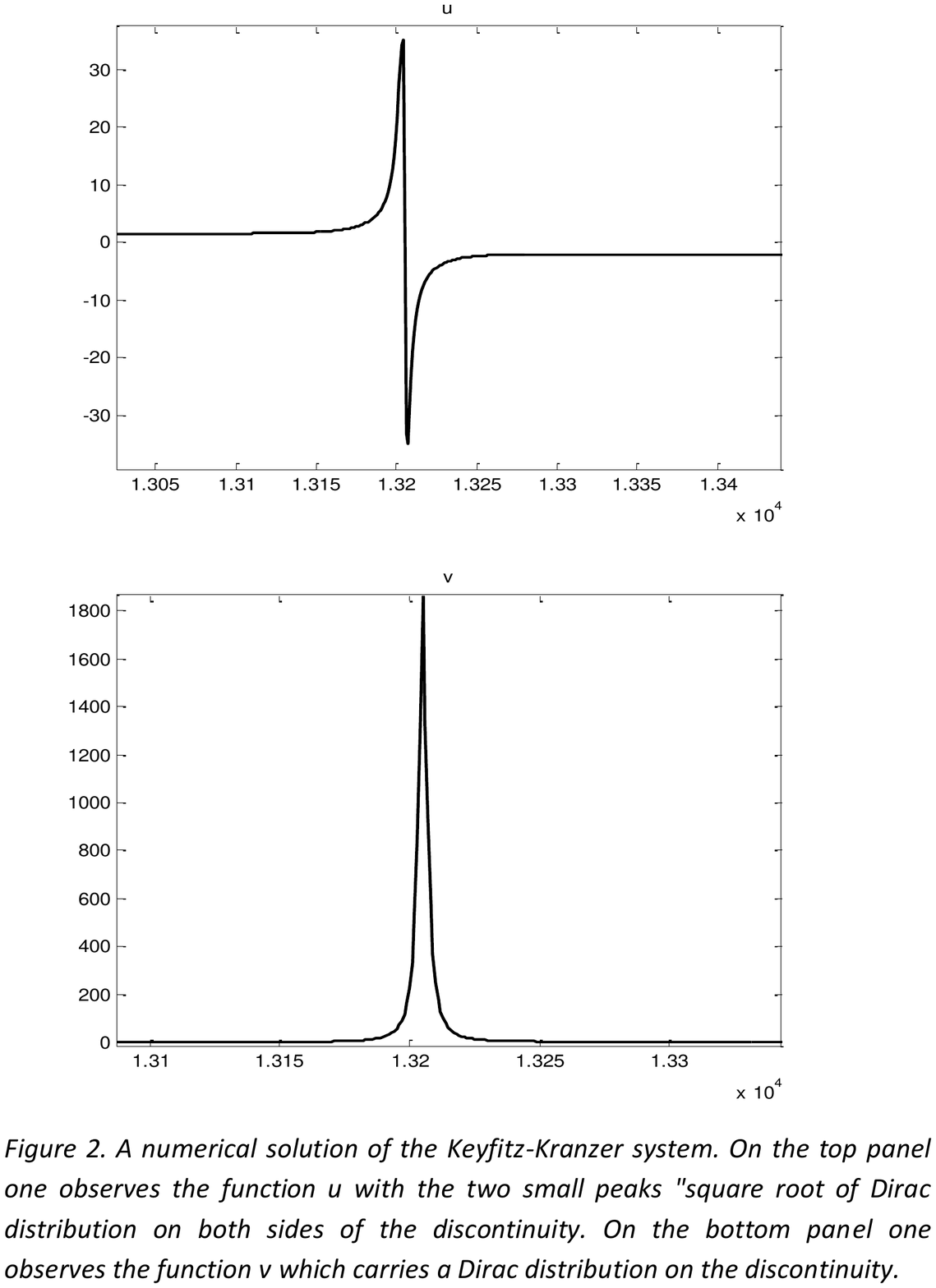}
\end{figure}
\\
where $\rho\in \mathcal{C}_c^\infty(]-1,+1[), \ \int\rho^2=1, \ \int\rho^3=0, \ x\rho(x)\leq 0, \ h$ a $\mathcal{C}^\infty$ monotone approximation of the Heaviside function, $h'\in \mathcal{C}_c^\infty(]-1,+1[), p >0$ a real number. The function $x\longmapsto \frac{1}{t\epsilon}\rho^2(\frac{\xi-s}{\epsilon})$ is a Dirac delta measure centered at $x=st$ (indeed $\int\frac{1}{\epsilon}\rho^2(\frac{\frac{x}{t}-s}{\epsilon})dx=t\int\rho^2(y)dy=t)$. Therefore the term $\frac{a}{\sqrt{\epsilon}} \rho(\frac{\xi-s}{\epsilon})$ has the form $a\sqrt{t} \sqrt{\delta}$.  Contrarily to the standard notation $\sqrt {\delta}$ this term has a negative and a positive small peak due to the fact $\int\rho^3(x)dx=0$ which forces the function $\rho$ to have positive and negative values that can be observed in figure 2.\\

As indicated above the function $u$ in figure 2 shows small peaks in form of "$\sqrt\delta$". At the limit the space step tends to 0 they tend to 0 in the sense of distributions. If we neglect them and state the solution as 
\begin{equation} u(x,t)=u_l+[u]H(x-ct),\end{equation}
\begin{equation}v(x,t)=v_l+[v]H(x-ct)+a(t)\delta(x-ct), \end{equation}
with  $a(t)$ a nonzero function of $t$, then in (6) $u_t$ shows a Dirac delta distribution from (8), $u^2-v$ shows a Dirac delta distribution since $a(t)\not=0$ therefore $(u^2-v)_x$ shows a term $a(t)\delta'$. Therefore $u_t+(u^2-v)_x$ shows terms $a(t)\delta' + b\delta$ with $a(t)\not=0$ which is impossible from (6) since $\delta$ and $\delta'$ are linearly independent: equation (6) cannot be satisfied. Similarly concerning (7):  $v_t$ shows a nonzero $\delta'$ term (the observed value $c$ is  different from 0) while from (8) $(\frac{u^3}{3}-u)_x$ cannot show  a $\delta'$ term:  (7) cannot be satisfied. \\

This proves that the $\sqrt\delta$ small peaks in $u$ cannot be omitted to explain that the observed $u,v$ are solutions of the equations (6, 7). In their presence the above contradictions disappear at once: for instance in (6) $u^2$ presents now a $\delta$ that can cancel the $\delta$ from $v$ and thus (6) can be satisfied.\\

\textbf{Conclusion. New mathematical objects such as 
different Heaviside functions
and square roots of Dirac delta distributions
should be introduced in mathematics in order to model real world situations.}\\

 Soon after 1950 when L. Schwartz "proved" that "multiplication of distributions is impossible in any mathematical context, possibly different of distribution theory" mathematical physicists believed that the problems with formal mathematical objects introduced by physicists that were not distributions should be solved  by modifying the formulation of physics so as to model physics completely within the setting of distributions. Today this viewpoint has completely failed. The unique option now is to introduce these objects imposed by the real world  into mathematics.  Such introduction was done about 1980 inside the research group of Professor Leopoldo Nachbin in Brazil.\\

\textbf{2. Introduction and study of "new generalized functions" inside the Leopoldo Nachbin group} (J. Aragona, H. A. Biagioni, J. F. Colombeau, Jo\~ao Toledo, R. Soraggi). About 1980 the group of Professor Leopoldo  Nachbin was studying holomorphic functions defined on locally convex spaces (and also approximation theory not considered here) \cite{Barroso1, Barroso2, Machado, Mujica, Zapata}.\\

Let $\Omega$ be an open set in $\mathbb{R}^n$, let $ \mathcal{D}=\mathcal{D}(\Omega)=\mathcal{C}_c^\infty(\Omega)$ denote the space of all $\mathcal{C}^\infty$ functions on $\Omega$ with compact support and $ \mathcal{E'}=\mathcal{E}'(\Omega)$ denote the space of all distributions on $\Omega$ with compact support. These spaces are equipped with their classical locally convex topologies defined in Schwartz distribution theory (the reader who does not know topological vector spaces can follow the sequel).  One can use as well holomorphic functions over  $\mathcal{E}'(\Omega)$ and $\mathcal{D}(\Omega)$. I chose $\mathcal{C}^\infty$ functions \cite{Colombeau0} because this is more general and will  give more numerous new mathematical objects. I did two very simple remarks.\\

Remark 1.  $\mathcal{C}^\infty(\mathcal{E'}) \subset\mathcal{C}^\infty(\mathcal{D})$.\\ 

Proof. This inclusion is obtained by the restriction map of a function defined on $\mathcal{E'}$ to a function defined on 
$\mathcal{D}$ since it is well known in Schwartz distribution theory that $\mathcal{D}$ is contained and dense in $\mathcal{E'}.\Box$\\

Remark 2. If we define an equivalence relation on $\mathcal{C}^\infty(\mathcal{E'}) $ by\\
$$\Phi_1 \mathcal{R} \Phi_2 \ \  iff  \  \ \Phi_1(\delta_x)= \Phi_2(\delta_x) \ \forall x\in \Omega$$
where $\delta_x$ denotes the Dirac measure at the point $x$, then the set of all equivalence classes is algebraically isomorphic to the classical algebra of all $\mathcal{C}^\infty$ functions on $\Omega$.\\

Proof. It is easy to prove that the map $x\longmapsto \delta_x$ is $\mathcal{C}^\infty$ from $\Omega$ into $\mathcal{E'}$ therefore the map $x\longmapsto \Phi(\delta_x)$ is $\mathcal{C}^\infty$ from $\Omega$ into $\mathbb{C}$. Now the surjection follows from the fact that if $f\in \mathcal{C}^\infty(\Omega)$, $f$ can be considered as the class of the map $\Phi$ defined by $\Phi(T)=<T,f> $ if $<,>$ is the duality brackets between 
$\mathcal{E'}$ and $\mathcal{C}^\infty(\Omega)$. The injection follows from the fact that if $f_1\not= f_2 \ \ \exists x\in \Omega \  /  \ f_1(x)\not= f_2(x).\Box$\\

Therefore from these two remarks one has the diagram\\

 \  \  \  \ \  \  \  \  \  \  \ \  \  \ $\mathcal{C}^\infty(\mathcal{E}'(\Omega)) \  \  \  \stackrel{quotient}{\longmapsto} \  \   \mathcal{C}^\infty(\Omega)$\\

  \  \  \  \ \  \  \  \  \  \  \ \  \  \ \  \ \ $\downarrow$ $restriction$\ \ \ \ \ \ \ \ \ \ \ \  \ \ \ \  \  \ \ \ \ \ \ \ \ \ \   
\  \  \\
  
 \  \  \  \ \  \  \ \  \  \  \ \  \  \ $\mathcal{C}^\infty(\mathcal{D}(\Omega)). \  \ \ \ \ \ \ \ \ \  \  \ \ \\$

The natural idea is to close the diagram.\\

\textbf{How to close the diagram?} If $\Phi\in\mathcal{C}^\infty(\mathcal{D}(\Omega))$, $\Phi(\delta_x)$ does not make sense. One naturally considers a sequence $\rho_{\epsilon,x} \rightarrow \delta_x,  \ \rho_{\epsilon,x}\in \mathcal{D}(\Omega), \ \int_{\mathbb{R}^n}\rho(x)dx=1, \ \rho_{\epsilon,x}(\lambda)=\frac{1}{\epsilon^n}\rho(\frac{\lambda-x}{\epsilon})$. Then $\rho_{\epsilon,x}\rightarrow \delta_x$  in $\mathcal{E}'(\Omega)$. A natural idea to extend the relation 
  $\mathcal{R}$ is to replace $\Phi(\delta_x), \ \Phi \in\mathcal{C}^\infty( \mathcal{E}'(\Omega))$, by $\Phi(\rho_{\epsilon,x}), \ \Phi \in \mathcal{C}^\infty(\mathcal{D}(\Omega))$. The statement "$\Phi(\rho_{\epsilon,x})\rightarrow 0$" would not lead to an algebra, so a more subtle definition had to be chosen. This was done  as follows in 1980-82 inside the group of Professor  Leopoldo Nachin. 
If $q\in \mathbb{N}$ set
\begin{equation}\mathcal{A}_q:=\{\rho \in \mathcal{D}(\mathbb{R}^n) \ / \int_{\mathbb{R}^n}\rho(x)dx=1,\int_{\mathbb{R}^n}x^i\rho(x)dx=0 \ \forall i\in \mathbb{N}^n, 1\leq |i|\leq q\}.\end{equation}

Then if $\Phi\in \mathcal{C}^\infty(\mathcal{E}'(\Omega))$,  $\Phi(\delta_x)=0$  is  equivalent   to  
 \begin{equation}  \forall  \ \rho\in \mathcal{A}_q  \forall K\subset\subset \Omega \  \exists C>0, \eta>0 \ / \  \Phi(\rho_{\epsilon,x})
=O(\epsilon^{q+1}) \ .\end{equation}

A proof can be found in \cite{ColombeauJMAA} p. 99-100, \cite{Colombeau1} p. 58-60. This characterization is basic since it puts in evidence an ideal  in the subalgebra of the algebra $\mathcal{C}^\infty(\mathcal{D}(\Omega))$
of all  $\mathcal{C}^\infty$ functions on $\mathcal{D}(\Omega)$ made of those $\Phi$ such that 
\begin{equation}\forall \ K\subset\subset \Omega \ \  \exists N\in \mathbb{N} \ / \ sup_{x\in K} |\Phi(\rho_{\epsilon,x})|=O(\frac{1}{\epsilon^N}) \ \forall \ \rho\in \mathcal{A}_1.\end{equation}
This subalgebra was called the subalgebra of the "moderate" elements of $\mathcal{C}^\infty(\mathcal{D}(\Omega))$, from the bounds in $O(\frac{1}{\epsilon^N})$.
 Further one has to insert all partial $x$-derivatives in these definitions. Finally it was possible to close the diagram:\\

 \  \  \  \ \  \  \  \  \  \  \ \  \  \ \  \ \ \ $\mathcal{C}^\infty(\mathcal{E}'(\Omega)) \  \ \ \ \ \  \  \stackrel{quotient}{\longmapsto} \ \ \ \ \ \  \   \mathcal{C}^\infty(\Omega)$\\

  \  \  \  \ \  \  \  \  \  \  \ \  \  \ \ \ \ \  \  \ \ $\downarrow$ \ \ \ \ \ \ \ \ \ \ \ \  \ \ \ \ \ \  \ \ \ \  \  \ \ \ \ \ \ \ \ \ \ \ \  \ \  $\downarrow $  
\  \  \\
  
 \  \    $suitable \ subalgebra \ of \  \mathcal{C}^\infty(\mathcal{D}(\Omega)) \  \stackrel{quotient}{\longmapsto}\ \ \mathcal{G}(\Omega)$.\\
 
 This new space $\mathcal{G}(\Omega)$ is a differential algebra i.e. presence of partial derivatives and of a multiplication with all the usual algebraic properties. We have the sequence of algebraic inclusions:\\
\begin{equation}\mathcal{C}^\infty(\Omega)\subset\mathcal{C}^0(\Omega)\subset\mathcal{D}'(\Omega)\subset\mathcal{G}(\Omega).\end{equation}
\  \  \\
The partial derivatives in $\mathcal{G}(\Omega)$ induce the partial derivatives in $\mathcal{D}'(\Omega)$ which is a vector subspace of $\mathcal{G}(\Omega)$. $\mathcal{C}^\infty(\Omega)$ is a faithful subalgebra of $\mathcal{G}(\Omega)$.\\

Two delicate points, that will be basic for the understanding of nonlinear generalized functions appear at once: what is the real situation of the algebra $\mathcal{C}^0(\Omega)$ in (13) and where are the different Heaviside functions we need to represent the real world? Up to now one can check that $\mathcal{C}^0(\Omega)$ is not a subalgebra of $\mathcal{G}(\Omega)$ and that powers in $\mathcal{G}(\Omega)$ of the Heaviside distribution $\in \mathcal{D}'(\Omega) \subset \mathcal{G}(\Omega)$ are not in  $\mathcal{D}'(\Omega)$. These two points are due to the Schwartz impossibility result we will consider later and will be (quite satisfactorily) solved below.\\

The detailed definitions that lead to the definition of $\mathcal{G}(\Omega)$ are as follows. First we recall the notation $D^k=\frac{\partial^{|k|}}{\partial x_1^{k_1}\dots\partial x_n^{k_n} }$.\\

The subalgebra of the "moderate" elements of $\mathcal{C}^\infty(\mathcal{D}(\Omega))$ is defined as:
$$\{\Phi\in\mathcal{C}^\infty(\mathcal{D}) \ /  \ \forall K\subset\subset\Omega \ \  \forall  k\in \mathbb{N}^n \ \exists N\in \mathbb{N} \ / \  \forall \rho\in \mathcal{A}_1 \ \exists C,\eta>0 \ / \ \ $$
$$ \ |D^k\Phi(\rho_{\epsilon,x})| \leq \frac{C}{\epsilon^N} \ \forall \epsilon\in ]0,\eta[, \forall x\in K\}.$$

The kernel defining the equivalence relation is defined as:
$$\{\Phi\in\mathcal{C}^\infty(\mathcal{D}) \ /  \ \forall K\subset\subset\Omega \ \  \forall  k\in \mathbb{N}^n \ \exists N\in \mathbb{N} \ / \forall q>N, \ 
 \forall \rho\in \mathcal{A}_q \ \exists C,\eta>0 \ / \ \ $$
$$ \ |D^k\Phi(\rho_{\epsilon,x})| \leq C\epsilon^{q-N} \ \forall \epsilon\in ]0,\eta[, \forall x\in K\}.$$
Then:
$$ \mathcal{G}(\Omega)=\frac{subalgebra \ of \ moderate \ elements}{above \ kernel}.$$

$ \mathcal{G}(\Omega) $ is often called the "full algebra" of nonlinear generalized functions.
 M. Grosser (Vienna) discovered later that the definition of the  kernel could be somewhat simplified but anyway these definitions are complicated with too many quantifiers. A subalgebra of $\mathcal{G}(\Omega)$, denoted by $\mathcal{G}_s(\Omega)$ and called "simplified" or "special", was constructed for this purpose. The canonical character of the inclusion of $\mathcal{D}'(\Omega)$ into $ \mathcal{G}(\Omega) $ is lost in the simplified case. Also, the simplified, and even the full case, can be formulated (in a very closely related form, not exactly equivalent) without the use of $\mathcal{C}^\infty$ or holomorphic  functions defined over locally convex spaces, which makes a  departure from the main theme of Prof. Leopoldo Nachbin group, but, later,  one came back to  $\mathcal{C}^\infty$ or holomorphic  functions defined over locally convex spaces for nonlinear generalized functions on manifolds. Similarly the classical theory of distributions can be presented with locally convex spaces or also even without infinite dimensional spaces. The simplified or special algebra $\mathcal{G}_s(\Omega)$ has been more used than the full algebra. It is defined by a quotient like the full algebra.\\

Reservoir of representatives=
$\{f:\Omega\times]0,1]\longmapsto\mathbb{C}, \ \mathcal{C}^\infty \ in \ x \ for \ any \ fixed \ \epsilon,  $ $$ \ / \  \forall K\subset\subset\Omega \ \  \forall  k\in \mathbb{N}^n \ \exists N,C>0 \ / \ |D^k f(x,\epsilon)| \leq \frac{C}{\epsilon^N} \ \forall \epsilon\in ]0,1[, \forall x\in K\}.$$
Kernel=
$$\{f\in Reservoir \ / \ \  \forall K\subset\subset\Omega \ \  \forall  k\in \mathbb{N}^n  \  \  \forall q \ large \ enough , \ 
  \ \exists C_q>0 \ / \ $$
$$ \ |D^k f(x,\epsilon)| \leq C_q\epsilon^{q} \ \forall \epsilon\in ]0,1[, \forall x\in K\}.$$
Then:
$$ \mathcal{G}_s(\Omega)=\frac{Reservoir}{Kernel}.$$

The full algebra as defined above is based on infinite dimensional analysis (holomorphic or $\mathcal{C}^\infty$ functions defined over locally convex spaces) in the same way as Schwartz presentation of distributions uses also infinite dimensional analysis (continuous linear functions defined over  locally convex spaces). The full algebra appears as a faithful extension of Schwartz presentation: a large functional space $\mathcal{C}^\infty(\mathcal{D}(\Omega))$ of nonlinear functions and a quotient that reduces on the small subspace $\mathcal{L}(\mathcal{D}(\Omega), \mathbb{C})$ of linear functions (= the distributions) to absence of quotient, and thus gives back exactly the definition of distributions as a particular case.  The simplified algebra does not use infinite dimensional analysis and is a direct extension of the sequential definition of distributions that we recall now.  To define the distributions on an open set $\Omega\subset \mathbb{R}^n$ following the sequential approach one first considers the set of all families $\{f_\epsilon \}_{\epsilon\in ]0,1]}$ of $\mathcal{C}^\infty$ functions on $\Omega$ such that 
$$ \forall \phi\in \mathcal{C}_c^\infty (\Omega) \ \ \int f_\epsilon(x)\phi(x)dx \rightarrow \ a \ limit \ in \ \mathbb{R} \  or \ \mathbb{C} \ when \ \epsilon \rightarrow 0.$$ 
This set forms obviously a vector space. On this vector space one defines the equivalence relation 
$$ \{f_\epsilon\} R \{g_\epsilon\} \ \Leftrightarrow \ \forall \phi \in \mathcal{C}_c^\infty(\Omega) \ \int (f_\epsilon-g_\epsilon)(x)\phi(x)dx\rightarrow 0 \ when \ \epsilon\rightarrow 0.$$
The distributions are the equivalence classes. These distributions are the same objects as the distributions defined by Schwartz using the theory of locally convex spaces. Of these two presentations one could ask: which is the best? The question has already no answer in distribution theory. The success of Schwartz has occulted the sequential presentation (Mikusinski and co-authors) but now one tends to drop the theory of locally convex spaces, thus getting closer and closer to the sequential approach.\\

Prof. Leopoldo Nachbin  made me publish 3 books in his collection "Notas de Matematica" published by North Holland Pub. Company, now Elsevier:\\
"Differential Calculus and Holomorphy ", 1982, to prepare the background of $\mathcal{C}^\infty$ functions over locally convex spaces,\\
"New Generalized Functions and Multiplication of Distributions", 1984, to introduce the full algebra $ \mathcal{G}(\Omega)$,\\
"Elementary Introduction to New Generalized Functions", 1985, to introduce the simplified algebra $ \mathcal{G}_s(\Omega)$. \\
 Prof. Leopoldo Nachbin presented me to Prof. L. Schwartz who supported this theory \cite{ColombeauNotes},  and  4 members of his group, Jorge Aragona, Hebe Biagioni, Jo\~ao Toledo and Roberto Soraggi   started working on the field of nonlinear generalized functions thus created \cite{AraBia,Bialivre}.\\

In his famous book "Th\'eorie des distributions", L. Schwartz claims: " multiplication of distributions is impossible in any mathematical context, possibly different from distribution theory". Since $\mathcal{D}'(\Omega)\subset \mathcal{G}(\Omega)$ canonically, or $\mathcal{D}'(\Omega)\subset \mathcal{G}_s(\Omega)$ non canonically, the algebraic structure of  $\mathcal{G}(\Omega)$ or $\mathcal{G}_s(\Omega)$  points out  a multiplication of distributions. Where is the mistake?
The first obvious answer is that the above claim is not a theorem, but the interpretation of a theorem. So one has to go back to the Schwartz theorem that we recall now.\\

Theorem (Schwartz 1954 \cite{Schwartz}). Let $A$ be an associative algebra containing the algebra $\mathcal{C}(\mathbb{R})$ as a subalgebra. Let us assume \\
\\
$\bullet $ the map $(x\longmapsto 1) \in \mathcal{C}(\mathbb{R})$ is unit element in $A$,\\
$\bullet$  there exists a map $D:A\longmapsto A$ called derivation such that\\
 $D(uv)=Du. v + u.Dv$, and
 $D|_{\mathcal{C}^1(\mathbb{R})}$ is the usual derivative.\\
\\
Then $D^2(x\longmapsto|x|)=0$ in $A$ instead of $2\delta$ as it should be.\\

One checks at once that all assumed properties hold in the context of  $\mathcal{G}(\mathbb{R})$, or $\mathcal{G}_s(\mathbb{R})$,  except that the classical algebra $\mathcal{C}(\mathbb{R})$ is not a subalgebra of $\mathcal{G}(\mathbb{R})$ or $\mathcal{G}_s(\mathbb{R})$,  although $\mathcal{C}^\infty(\mathbb{R})$ is. Let $f,g\in \mathcal{C}(\Omega)$. Then the new product $f\bullet g\in \mathcal{G}(\Omega)$ is the class, modulo the kernel defining $\mathcal{G}$,  of the map $\mathcal{D}(\Omega)\longmapsto \mathbb{C}$ whose restriction to the set of the $\rho_{\epsilon,x}$ is
$$\rho_{\epsilon,x} \longmapsto \int f(x+\epsilon \mu)\rho(\mu)d\mu. \int g(x+\epsilon \mu)\rho(\mu)d\mu.$$

On the other hand the classical product $f.g$, when considered in $\mathcal{G}(\Omega)$ through the inclusion map $\mathcal{C}(\Omega) \subset \mathcal{G}(\Omega)$, is the class of 
$$\rho_{\epsilon,x} \longmapsto \int f(x+\epsilon \mu)g(x+\epsilon \mu)\rho(\mu)d\mu.$$
The difference  $(f\bullet g -f.g)(\rho_{\epsilon,x})\rightarrow 0$ when $\epsilon\rightarrow 0$. From the two definitions of the kernel above (in the full or simplified case) this is not a sufficient smallness to ensure that the difference $(f\bullet g -f.g)$ is in any of the kernels since a fast tendance to $0$ (of the form $o(\epsilon^q) \ \forall q$) is requested. This why the two products are different.\\

Why  $f\bullet g = f.g$ when $f,g\in \mathcal{C}^\infty(\Omega) $? This follows from Taylor's formula and the definition of the sets $\mathcal{A}_q$. Indeed\\
$\int[f(x)+\epsilon \mu f'(x)+\dots+\epsilon^q\mu^qf^{(q)}(x)+\epsilon^{q+1}\mu^{q+1}r(x,\epsilon)]\rho(\mu)d\mu = f(x) \int\rho(\mu)d\mu+\epsilon f'(x)\int\mu \rho(\mu)d\mu+\dots+\epsilon^q f^{(q)}(x)\int \mu^q\rho(\mu)d\mu+\epsilon^{q+1}\int r(x,\epsilon) \mu^{q+1}\rho(\mu)d\mu= f(x)+\epsilon^{q+1}\int r(x,\epsilon) \mu^{q+1}\rho(\mu)d\mu=f(x)+O(\epsilon^{q+1})$ if $\rho\in \mathcal{A}_q$.\\

Therefore if $f,g\in \mathcal{C}^\infty$ and $\rho\in \mathcal{A}_q$
$$(f\bullet g-f.g)(\rho_{\epsilon,x})=O(\epsilon^{q+1})$$
which shows that the map (rather, the map  defined on $\mathcal{D}(\Omega)$ whose restriction to the set of all $\rho_{\epsilon,x}$ is):
$$\rho_{\epsilon,x}\longmapsto (f\bullet g-f.g)(\rho_{\epsilon,x})$$
is in the kernel defining $ \mathcal{G}(\Omega)$ if $f,g$ are $\mathcal{C}^\infty$ functions.\\

If $f,g \in \mathcal{C}^n(\Omega)$ one has only
$$(f\bullet g-f.g)(\rho_{\epsilon,x})=O(\epsilon^{n+1})$$
if $\rho\in \mathcal{A}_q$, whatever $q>n$, therefore $f\bullet g-f.g\not\in$ kernel.\\

This defines various concepts of smallness in  $\mathcal{G}(\Omega), \  \mathcal{G}_s(\Omega)$. A convenient concept which generalizes  the concept defining the distributions is called association and is defined as follows:\\

\textit{Definition. We say that $G_1,G_2\in \mathcal{G}(\Omega)$ are associated iff they have respective  representatives $G_i(\rho_{\epsilon,x}), i=1,2$  (then it works for arbitrary representatives) such that} 
\begin{equation} \forall \psi \in \mathcal{D}(\Omega), \ \forall \rho \in \mathcal{A}_N, N \ large\ enough,\ \int_{\Omega}(G_1(\rho_{\epsilon,x})-G_2(\rho_{\epsilon,x}))\psi(x)dx \rightarrow 0 \end{equation}
\textit{when} $\epsilon\rightarrow 0$.\\

We note $$G_1\approx G_2$$
to express that $G_1$ and $G_2$ are associated. We have
\begin{equation} H^n\approx H \ \ \forall n, \ \ \sqrt{\delta} \approx 0,\end{equation} if $H$ is the Heaviside function and $\delta$ the Dirac delta distribution. From the definition of association 
\begin{equation} G_1\approx G_2  \ \ implies \ \ \frac{\partial}{\partial x}G_1\approx \frac{\partial}{\partial x}G_2.\end{equation}
But
\begin{equation} G_1\approx G_2  \ \ does \ not \ imply \ \  G.G_1=G.G_2:\end{equation}
it suffices to consider $G_1=G=\sqrt{\delta}, \  G_2=0$.\\

Of course for every distribution $T$ there are many generalized functions $G\in \mathcal{G}(\Omega)$ associated to $T$. In the sense of the definition of association they have the same "macroscopic aspect" as $T$, but can differ completely from $T$ by deeper properties: $\sqrt{\delta}$ and 0 for instance.\\

 \textbf{Conclusion. In $\mathcal{G}(\Omega)$ we have two very different concepts that play the role of the classical  equality:}\\

$\bullet$ \textbf{of course the "true" equality in $\mathcal{G}(\Omega)$ which is coherent with all operations: multiplication, differentiation,}\\

$\bullet$ \textbf{the association, which is not an equality since different generalized functions can be associated, but extends faithfully the concept of equality of distributions to this more subtle setting. The association is coherent with differentiation but incoherent with multiplication, as a reproduction of the properties of the equality of distributions.}\\

Note a property that could look as a paradox: we have an inclusion of the space of distributions into the differential algebras of nonlinear generalized functions (canonical in the full case, noncanonical in the simplified case); then two generalized functions which are in the image of the space of distributions are associated iff they are  equal in  $ \mathcal{G}(\Omega)$: this follows from the fact that the other associated generalized functions are excluded since they are not in  this image.\\

One can interpret the classical products of distributions through the association:\\

$\bullet$ the Schwartz product $f.T, \ f\in \mathcal{C}^\infty(\Omega), T\in \mathcal{D}'(\Omega)$ and the new product $f\bullet T \in \mathcal{G}(\Omega)$  are associated \cite{Colombeau1,ColombeauJMAA};\\

$\bullet$ the same holds for the Mikusinski product \cite{Colombeau1}, and more particularly for the Hormander product \cite{Colombeau1,Tysk} due to a result  of \cite{Tysk}.\\

$\bullet$ More generally the following coherence between the classical calculations and the calculations in 
$\mathcal{G}(\Omega)$ has always been observed: \textit{whenever a sequence of calculations makes fully sense in classical mathematics or in distribution theory, the analog calculations in the $\mathcal{G}$-context give a result associated to the classical result.}\\

Now one can understand the Schwartz impossibility result: the assumption that the classical algebra $\mathcal{C}(\mathbb{R})$ is a subalgebra of the algebra $A$ is stronger than needed to ensure coherence of the new calculus in $A$ and the classical calculus. Indeed coherence concerns only the calculations that fully make sense, from their beginning to their end, in classical mathematics and distribution theory. In the case $A=\mathcal{G}(\mathbb{R})$ this coherence is obtained: with the new product of continuous functions one obtains in $\mathcal{G}(\mathbb{R})$ a result which is simply more precise than the classical result and it suffices to drop this refinement (using the association) to have exactly the classical result. The Schwartz assumption that $\mathcal{C}(\mathbb{R})$ would be a subalgebra of the algebra $A$ means that one could use the classical product even outside the domain of classical calculations. For instance, in a sequence of calculations such that one at least does not make sense classically or in distribution theory, consider a previous calculation of the sequence that makes sense classically such as a product of two continuous functions: if the Schwartz assumption $\mathcal{C}(\mathbb{R})$ subalgebra of $A$ would be satisfied one could use the classical product: this is not the case  because of the  nonclassical calculations that follow this classical calculation. One should use in this case the new product of continuous functions: an infinitesimal difference such as $H^2-H$ can be multiplied by an infinite quantity such as $H'$ and give a significant result. Since this difficulty concerns a sequence of calculations in which one at least does not make sense within classical mathematics or  distribution theory the use of classical calculations in a part of them is more than needed for coherence between classical mathematics and the calculations in the new setting. In short in the context of $\mathcal{G}(\Omega)$ one has complete coherence but only in form of the minimum acceptable to have coherence,  and  the Schwartz assumption that $\mathcal{C}(\mathbb{R})$ would be a subalgebra of $A$ would give more than this minimum. It is precisely this subtle fact that  makes the difference between possibility and impossibility.\\

The association also serves to interpret weak solutions of PDEs:\\

$\bullet$ Let us recall the definition of a weak asymptotic method as stated by Danilov-Omel'yanov and Shelkovich \cite {DOS} as an extension of Maslov asymptotic methods: for simplicity consider the equation $u_t+f(u)_x"="0$ where $"="$ means that the sense of this equality is not yet fixed. Then a weak asymptotic method is a sequence $(u^\epsilon)_{\epsilon \rightarrow 0}$ such that $\forall \psi \in \mathcal{D}(\mathbb{R}) \  \forall t \  \int ((u^\epsilon)_t\psi -f(u^\epsilon)\psi_x)dx\rightarrow 0$ as $\epsilon\rightarrow 0$. One recognizes the similarity with the definition of the association in $\mathcal{G}_s(\Omega)$ in variable $x$ for fixed $t$:  $\forall \psi \in \mathcal{D}(\mathbb{R}) \  \forall t \  \int ((u^\epsilon)_t\psi +(f(u^\epsilon))_x\psi)dx\rightarrow 0$ as $\epsilon\rightarrow 0$. The difference lies in that in the weak asymptotic method $f(u^\epsilon)$ needs not be $\mathcal{C}^\infty$.
 Indeed, provided the family $u^\epsilon$ is made of $\mathcal{C}^\infty$ functions and has bounds in $\frac{const}{\epsilon^N} $  on compact sets of $\Omega$ and for all partial derivatives, the two concepts are same: starting from any weak asymptotic method these two requirements (i.e. $ \mathcal{C}^\infty$  and bounds) can be obtained by convolution and by reindexation in $\epsilon$ respectively. Therefore \\

\textbf{There exists a weak asymptotic method iff there exists a solution in the sense of association.}\\

 The concept of weak asymptotic method as defined in \cite{DOS} has been at the origin of many important investigations on solutions of systems of conservation laws in recent years \cite{Albeverio, AlbeShelk, Danilov, DOS, Mitrovic, Mitrovic2,DanilovO1, DanilovO2, DanilovSh, Shelkovich2, Shelkovich3, Mitrovic3, Omel'yanov, Panov, ShelkovichRMS, Shelkovichmat, Shelkovich1} \dots.\\

$\bullet$ In presence of $L^p$ bounds, $1\leq p \leq \infty,$, anyway  needed to produce a measure valued solution,  it has been proved that\\

 \textbf{The existence of a solution in the association sense is equivalent to existence of a measure valued solution}, \cite{Alecsic}.\\

The need of solutions in the association sense is clearly shown by the following well known remark: the equation $u_t+uu_x"="0$ has not same jump conditions as the equations obtained by multiplying it by $u$, for instance $uu_t +u^2u_x"="0$. The interpretation is that the equation stated in $\mathcal{G}$ in the strong form $u_t+uu_x=0$ (i.e. with the "true" equality in $\mathcal{G}$) has no shock wave solution since multiplication by $u$ would give same solution; the equation stated in $\mathcal{G}$ in its weak form $u_t+uu_x\approx 0$ (i.e. with the association)
has shock wave solutions but multiplication by $u$ is forbidden, which is put in evidence by the fact it changes the solutions.\\

\textbf{Conclusion. In order to multiply distributions and to create new needed mathematical objects: $H^n\not = H^p$  if $n\not= p$, $\sqrt{\delta}\not= 0$, a context of nonlinear generalized functions has been introduced inside the Leopoldo Nachbin group (1980-1986): Jorge Aragona and Jo\~ao Toledo (USP), Hebe  Azevedo Biagioni (UNICAMP), Jean Fran\c cois Colombeau (visitor at UNICAMP, USP, UFRJ), Roberto Soraggi (UFRJ). One has defined and studied two differential algebras (and variants)\\
$$\mathcal{G}_s(\Omega)\subset \mathcal{G}(\Omega)$$
based on a more refined concept of equality (asymptotic behavior in $o(\epsilon^q) \forall q$) instead of mere tendance to 0. Dropping the refinement of this asymptotic behavior in favor of a tendance to 0 gives the association (plus use of test functions in the standard definition of association). The concept of  association  had already been  defined earlier in various circumstances  by numerous authors such as Mikusinski (products of distributions), Lojasiewicz (point values of distributions), Maslov, Danilov, Omel'yanov, Shelkovich (asymptotic methods).}\\

\textbf{3. Direct subsequent applications after 1986}. After 1986 the theory evolved outside the Leopoldo Nachbin group; in particular let us quote:\\
$\bullet$ in France (J.F. Colombeau, A.Y. Le Roux and coauthors) under supervision of the technical branch of the Ministry of Defense (DRET): simulation of strong collisions for armor, wave diffusion for furtivity.\\
$\bullet$ in Austria (M. Oberguggenberger and coauthors  M. Kunzinger, C. Garetto, G. Hormann,\dots) and Serbia (S. Pilipovic and coauthors): solutions of linear and nonlinear Cauchy problems with = in $\mathcal{G}$, propagation of singularities for linear PDEs with irregular coefficients (numerous articles).\\
$\bullet$ in Austria (M. Grosser, M. Kunzinger, R. Steinbauer) and England (J. Vickers): generalized functions on manifolds and General Relativity, see the research-expository paper \cite{Nigsch} in this volume.\\
$\bullet$ in Brazil (J. Aragona, S.O. Juriaans and coauthors): algebraic properties and use of the Scarpalezos sharp topology \cite{AragonaFernandez, AragonaFerJu1, Novo, AragonaFerJu2, AragonaFerJu3} and also \cite{Garetto, Garetto2,GarettoV, Vernaeve, Scarpatopo, Scarpatopo2}.\\
$\bullet$ Stochastic equations (F. Russo, M. Oberguggenberger, D. Rajter, S. Pilipovic, D. Selesi, P. Catuogno, C. Olivera,\dots) \cite{OberRusso, Russo, MPS, CaO1,CaO2, RS}.\\ 
$\bullet$  Connections with nonstandard analysis have been developped by various authors \cite{Hoskins, Todorov, Vernaeve, Ober}, \dots, see  \cite{Todorov} in this volume.\\

The theory of nonlinear generalized functions is now presented in many books or research expository papers \cite{AraBia, Bialivre, Christiakov, Colombeau1,  Colombeau2,  Colombeau4,  ColombeauBAMS, Egorov, GKOSmemoirs, GKOS, Nedel, Ober, Rosinger}.
 Now we sketch some  applications that permit a better understanding of the theory and its use.\\

\textbf{Symmetric linear hyperbolic systems (M. Oberguggenberger).} One considers  linear symmetric hyperbolic systems

$$ \frac{\partial U}{\partial t}=\sum_{i=1}^{n}A_i(x_1,\dots,x_n,t)\frac{\partial U}{\partial x_i}(x_1,\dots,x_n,t)+B(x_1,\dots,x_n,t)U(x_1,\dots,x_n,t)+$$\begin{equation}C(x_1,\dots,x_n,t),\end{equation}
where the $A_i$ are real symmetric matrices with coefficients in $\mathcal{G}(\mathbb{R}^{n+1})$ which further have a specific growth (in form of a bound $log(\frac{1}{\epsilon})$)  when $\epsilon\rightarrow 0$, see \cite{Lafon} for  details.\\

The authors obtain existence and uniqueness of solutions of the Cauchy problem with = in the $\mathcal{G}$ context. Further the solution is associated to the classical solution when the later exists.\\

This result has permitted to consider the case of discontinuous coefficients, very important in the applications since it models wave propagation in an heterogeneous medium and it has been at the origin of numerous works. 
The proof consists in applying to representatives a classical proof of the $\mathcal{C}^\infty$ case that is refined  more in detail so as to prove that the regularized solutions have the bounds in $\frac{1}{\epsilon^N}$ and to prove uniqueness from  bounds in $o(\epsilon^q) \forall q$. This technique of working on representatives by refining a classical $\mathcal{C}^\infty$ proof and using the definition of nonlinear generalized functions as a quotient has been applied in many cases and it has permitted to extend many $\mathcal{C}^\infty$ results to the case of nonlinear generalized functions.\\

\textbf{A nonlinear very dissipative parabolic equation.}  Another example of the above technique consisting in refining classical $\mathcal{C}^\infty$ proofs has been used in \cite{ ColombeauLanglais, Langlais, AragonaGarcia} for the parabolic equation\\
\begin{equation} u_t-\Delta u+u^3=0,\end{equation}
\begin{equation} u(.,0)=f \in \mathcal{E}'(]-1,1[),\end{equation}
\begin{equation} u(x,t)=0 \ on \ \partial \Omega\times [0,+\infty[.\end{equation}
When the initial condition is a distribution with compact support in $]-1,+1[$ one obtains  existence and uniqueness of a solution with equality in the algebra $\mathcal{G}_s([-1,+1]\times [0,+\infty[)$ which is defined like the case of generalized functions over an open set (representatives of class $\mathcal{C}^\infty$ in the interior, continuous as well as all derivatives up to the boundary, and uniform bounds on compact subsets of $[-1,+1]\times[0,+\infty[$). This solution is associated to the classical solution when the later exists.\\

An interesting point is as follows:  it had been noticed \cite{Brezis} that if a Dirac delta distribution is approximated as usual by a sequence $(\delta^\epsilon)$ then the  classical $\mathcal{C}^\infty$ solutions $u^\epsilon$ corresponding to the initial conditions  $\delta^\epsilon$ tend to 0 uniformly on compact subsets of $[-1,+1]\times ]0,+\infty[$ when $\epsilon\rightarrow 0$. If as usual the solution $u$ is defined as a limit of the $u^\epsilon$ then $u=0$, and therefore the initial condition defined as limit of $u(.,t)$ when $t\rightarrow 0$ should be the null function, while we started with the Dirac delta distribution as initial condition. This paradox suggests  to consider that problem (19-21) with Dirac distribution as initial condition has no solution. In the $\mathcal{G}$ context the explanation of the paradox is clear: the solution $u\in  \mathcal{G}_s([-1,+1]\times [0,+\infty[)$ is associated to 0 in   $\Omega\times ]0,+\infty[$, but it is not equal to 0, and its initial condition which in the $\mathcal{G}$ context is defined as its restriction to $t=0$ is by construction the given distribution on  $]-1,+1[$.\\

\textbf{Weak solutions in the context of nonlinear generalized functions.} The solutions considered in the two above sections were strong solutions in the sense that they were solutions with the equality in the algebra of nonlinear generalized functions. Unfortunately many basic equations do not have such solutions and one is forced to consider also solutions in the association sense. This has already been noticed for the shock waves solutions  of the equation $u_t+uu_x"="0$ - where we recall that the symbol $"="$ means that the sense to be given to this equality is not yet fixed- that we are forced to interpret as $$u_t+uu_x\approx0.$$ This is equivalent to a weak asymptotic method, see above, with weak derivative both in $x$ and $t$, i.e.  a family $(u^\epsilon)_\epsilon$ of measurable functions in the variables $x,t$ such that, in the case of the  equation $u_t+(f(u))_x\approx 0$,
\begin{equation} \forall \psi \in \mathcal{C}_c^\infty(R^2) \ \ \int(u^\epsilon\psi_t+f(u^\epsilon)\psi_x)dxdt\rightarrow 0 \ \ when\ \  \epsilon\rightarrow 0.\end{equation}

A shock wave in the usual form of two constant values $u_l,u_r$ separated by a discontinuity moving with constant velocity $c$ is represented by the formula
\begin{equation} u(x,t)=u_l +[u] H(x-ct),\end{equation}
where $[u]:=u_r-u_l$ is the jump of $u$ and $H$ is any Heaviside function in $\mathcal{G}(\mathbb{R})$ or $\mathcal{G}_s(\mathbb{R})$. What is precisely a Heaviside generalized  function? \\ 

A Heaviside generalized function $H$ is defined as an element of $\mathcal{G}(\mathbb{R})$ which has a  representative $(x,\epsilon)\longmapsto H(x,\epsilon)$ uniformly bounded in $\epsilon$ and $x$ such that $H(x,\epsilon)$ tends to 0 if $x<0$ and 1 if $x>0$ in a sense that can 
depend on authors, but that implies that $\forall n \  H^n\approx H$. All Heaviside  functions are associated and $\int_{-\infty}^{+\infty}H'(x,\epsilon)dx=1$.\\

 It will be very important for the sequel to remark that $H$ should not be the privilieged Heaviside function in $\mathcal{D}'(\mathbb{R})$ from the  (canonical or not) injections of $\mathcal{D}'(\mathbb{R})$ into the two algebras $\mathcal{G}(\mathbb{R})$ or $\mathcal{G}_s(\mathbb{R})$ of nonlinear generalized functions. Indeed we shall deal with several physical variables such as density $\rho$, velocity $u$,\dots, that in general cannot  be represented by the same Heaviside function (recall the example of the elastoplastic shock waves) so that we will need different Heaviside functions and there is no reason to priviliege one physical variable. We will consider equations in divergence form also called conservative equations because they express conservation of some quantity, and equations not in divergence form also called nonconservative equations.\\

Inserting (23) into the equation $u_t+(f(u))_x\approx 0$  one obtains at once the classical jump condition
\begin{equation} c=\frac{f(u_r)-f(u_l)}{u_r-u_l}=\frac{[f(u)]}{[u]}.           \end{equation}
The detailed calculation is easy: inserting (23) into (24) one obtains $-c[u]H'+(f(u_l+[u]H))'\approx 0$. One can integrate an association as this is classically the case for the equality in distribution theory; one obtains: $ -c[u]H+f(u_l+[u]H)+const \approx 0.$ Using $H(-\infty)=0$ one obtains $const=-f(u_l)$; then $H(+\infty)=1$ gives the jump condition (24).\\

One obtains easily: \textit{$u$ in (23) is solution in the association sense iff the jump condition  (24) holds}. A solution in the association sense does not bring any information on the Heaviside function $H$: all Heaviside functions fit as well: although there is uniqueness of jump condition,  and so uniqueness,  modulo association, 
of solutions of form (23) for the Cauchy problem, there is no uniqueness in the sense of equality in the $\mathcal{G}$ context since the function $H$ is not fixed.\\

 Now we know that a nonlinear scalar equation  has no shock wave in the sense of equality in  $\mathcal{G}$ (same as above for $u_t+uu_x=0$), so let us turn to a system of two equations. We will distinguish systems in divergence (or conservative) form and systems in nondivergence (or nonconservative) form since the shock wave theory is quite different in each case. \\

Systems in divergence form are systems 
\begin{equation}u_t+(f(u,v))_x"="0, \end{equation}
\begin{equation}v_t+(g(u,v))_x"="0, \end{equation}
i.e. they have the very particular feature that they are sums of derivatives. Nonconservative systems are systems which are not in divergence form such as the system 
\begin{equation}u_t+(u^2)_x "=" 0, \end{equation}
\begin{equation}\sigma_t+u.\sigma_x "=" 0, \end{equation}
because of the term $u.\sigma_x$.
 The best known systems are in divergence form, such as the famous system of ideal gases that we state here in one dimension:

\begin {equation} \rho_t+(\rho u)_x"="0,\end{equation}
\begin {equation}(\rho u)_t+(\rho u^2)_x+ p_x"="0,\end{equation}
\begin {equation}(\rho e)_t+(\rho e u)_x+(pu)_x"="0,\end{equation}
\begin{equation}  p"="(\gamma-1)(\rho e -\rho \frac{u^2}{2}),\end{equation}
\\
 where $\rho, u, p, e$ denote respectively the density, the velocity, the pressure and the  total energy per unit mass; $\gamma$ is a constant =1.4 for the air in usual conditions. \\
 
For systems in divergence form stated with the association $$(u_i)_t+(f_i(u_1,\dots,u_n))_x\approx 0$$ we have noticed (23, 24) that shock waves  $u_i(x,t)=(u_i)_l+[u_i]H_i(x-ct), (u_i)_l, [u_i] \in \mathbb{R}, H_i$  Heaviside generalized functions, satisfy the jump conditions $c=\frac{[f]}{[u_i]}$. For systems in nondivergence form the situation is completely different: consider for instance the equation $u_t +uv_x-u_x\approx 0$ and set 
\begin {equation}u(x,t)=u_l+[u]H(x-ct), \ \ v(x,t)=v_l+[v]K(x-ct)\end{equation}
where $H,K$ are two Heaviside functions.
  Insertion of (33) into the equation gives $-c[u]H'(\xi)+u_l [v]K'(\xi) +[u][v]H(\xi)K'(\xi)-[u]H'(\xi)\approx 0$ as a generalized function of the variable $\xi=x-ct$. Next one would like to integrate in $\xi$: this is impossible because $\int H(\xi)K'(\xi) d\xi$ is unknown and can take a continuum of values according to $H$ and $K$. Therefore we do not obtain well defined jump conditions. This problem is important and will be treated in detail in the sequel.\\

\textbf{Examples of nonconservative systems.} We give three classical examples of nonconservative systems in which the origins of the  nonconservative equations  have  different natures: the state law for the elasticity system of large deformation \cite{ColombeauJMP}, the relativity equations \cite{Sachs} for the system of primordial cosmology \cite{Coles}, the statement of the basic laws of physics for the multifluid system \cite{Wendroff, Monk}. For convenience the systems are given in 2-D since the 3-D systems are obvious extensions of the 2-D systems.\\

\textit{Elasticity for large deformations.} The following nonconservative equations (38-40) stem from a differential form of Hooke's law. In elasticity the nonconservative terms in (38) are $ v.(s_{11})_y$ and $(u_y-v_x).s_{12}$. In elastoplasticity  $\mu$  is null in the plastic region, \cite{ColombeauJMP}, which introduces other "nonconservative terms" in the equations  (38-40). This system is used for the numerical simulations of strong collisions, in particular for the design of armor \cite{ColombeauJMP, Colombeau4, Bialivre}. \\

The model in 2-D is (stating = in place of "="):
\begin {equation} \rho_t+(\rho u)_x+(\rho v)_y=0,\end{equation}
\begin {equation}(\rho u)_t+(\rho u^2)_x+(\rho uv)_y+(p-s_{11})_x-(s_{12})_y=0,\end{equation}
\begin {equation}(\rho v)_t+(\rho uv)_x+(\rho v^2)_y+(p-s_{22})_y-(s_{12})_x=0,\end{equation}
\begin {equation}(\rho e)_t+(\rho eu)_x+(\rho ev)_y+((p-s_{11})u)_x+((p-s_{22})v)_y-(s_{12}v)_x-(s_{12}u)_y=0,\end{equation}
\begin {equation}(s_{11})_t+u.(s_{11})_x+v.(s_{11})_y-\frac{4}{3}\mu u_x-\frac{2}{3}\mu.v_y-(u_y-v_x).s_{12}=0,\end{equation}
\begin {equation}(s_{22})_t+u.(s_{22})_x+v.(s_{22})_y+\frac{2}{3}\mu v_x-\frac{4}{3}\mu.u_y+(u_y-v_x).s_{12}=0,\end{equation}
\begin {equation}(s_{12})_t+u.(s_{12})_x+v.(s_{12})_y-\mu v_x+\mu.u_y-\frac{1}{2}(u_y-v_x).(s_{11}-s_{22})=0,\end{equation}
\begin {equation} p=\Phi(\rho, e-\frac{u^2+v^2}{2}),\end{equation}
where $\rho=$ density; $(u,v)$=velocity vector, $p=$pressure, $e$=density of total energy; $(s_{ij})$=stress tensor ($s_{ij}=s_{ji}$), $\mu$= shear modulus which can be considered as constant in  elasticity, $\Phi$ = a state law.\\\\

\textit{Primordial cosmology.} According to the standard big bang theory, in its infancy the universe was "small" and the photons were very energetic (small wavelength). Therefore radiation dominated baryonic matter from well known quantum effects. Later the expansion of the universe increased the wavelength of the photons, and baryonic matter became free as it is today. The equations in the radiation dominated universe are very important since the seeds of galaxies were formed from black matter at this epoch, which provided gravitational wells for baryonic matter  when radiation domination ceased. The equations below are those in use by cosmologists \cite{Coles}. They are issued from the equations of relativity, for which the lack of divergence form  is a well known problem \cite{Sachs}: the Euler equations (43, 44) are not in divergence form. The equations are:\\

\begin {equation} \rho_t+((\rho+\frac{p}{c^2}) u)_x+((\rho+\frac{p}{c^2}) v)_y=0,\end{equation}
\begin {equation}(\rho+\frac{p}{c^2}) ( u_t+uu_x+uv_y)+ p_x+(\rho+\frac{p}{c^2})\Phi_x=0,\end{equation}
\begin {equation}(\rho+\frac{p}{c^2}) ( v_t+vu_x+vv_y)+ p_y+(\rho+\frac{p}{c^2})\Phi_y=0,\end{equation}
\begin {equation}\Delta \Phi=4\pi G(\rho+\frac{3p}{c^2})\end{equation}
\begin{equation}  p=\frac{1}{3}\rho c^2 ,\end{equation}
\\
 where $\rho, (u,v), p, $ denote respectively the density of energy, the velocity vector, the radiation pressure; $c$ is the velocity of light.\\

\textit{Multifluid flows.} The system below is called the "standard one pressure model". It models a mixture of two immiscible fluids which are strongly mixed, for instance numerous bubbles of one fluid into the other, for which there is only one pressure due to the mixture. The extremely complicated geometry of the surface separating the two fluids is replaced by an averaging from the respective volumic proportions $\alpha_1$ and $\alpha_2$ of the fluids, which introduces nonconservative terms. The model is derived and discussed in \cite{Wendroff, Monk}. It is used in offshore petroleum extraction to study the "gas kick", i.e. a destructive possible explosion when a mixture of oil and natural gas reaches the surface, see the extremely rich literature on the subject in \cite{Monk, Evje, EvjeFjelde,EvjeFlatten,EvjeFriis,EvjeCommunications, Toumi}\dots. This system is also used for the cooling of nuclear power stations (mixture of boiling water and vapor of water). Equations (49-52) are nonconservative.
 The 1-D model of a mixture in a tube 
(index 1 for fluid 1 and index 2 for fluid 2) is \cite{Wendroff, Monk,Toumi}:
\begin {equation}(\rho_1 \alpha_1)_t+(\rho_1\alpha_1 u_1)_x=0,\end{equation}
\begin {equation}(\rho_2 \alpha_2)_t+(\rho_2\alpha_2 u_2)_x=0,\end{equation}
\begin {equation}(\rho_1 \alpha_1u_1)_t+(\rho_1\alpha_1( u_1)^2)_x+\alpha_1p_x+\tau=\rho_1\alpha_1g,\end{equation}
\begin {equation}(\rho_2 \alpha_2u_2)_t+(\rho_2\alpha_2( u_2)^2)_x+\alpha_2p_x-\tau=\rho_2\alpha_2g,\end{equation}
\begin {equation}(\rho_1 \alpha_1e_1)_t+(\rho_1\alpha_1 u_1e_1)_x+(\alpha_1pu_1)_x+(p\alpha_1u_1)_x+p.(\alpha_1)_t=\rho_1\alpha_1u_1g,\end{equation}
\begin {equation}(\rho_2 \alpha_2e_2)_t+(\rho_2\alpha_2 u_2e_2)_x+(\alpha_2pu_2)_x+(p\alpha_2u_2)_x+p.(\alpha_2)_t=\rho_2\alpha_2u_2g,\end{equation}
where, if $i=1,2$, $\alpha_i$ is the volumic proportion of phase $i$ in the mixture, therefore, 
 \begin {equation}\alpha_1+\alpha_2=1,\end{equation}
 $\rho_i,u_i,e_i$ are the density, velocity and  total energy of phase $i$ per unit mass, $p$ is the  pressure  in the mixture, $\tau $ denotes the interface-momentum exchange term and $g$ is the component of gravitational attraction in the direction of the tube (which is not always vertical). The system is completed by the two state laws of the liquid and the gas: \begin {equation}p=F_i(\rho_i, e_i-\frac{u_i)^2}{2}), \ i=1,2.\end{equation}

\textbf{Jump conditions for nonconservative systems:} 
\textbf{the use of equality and association in the calculations of jump conditions for nonconservative systems.} Consider a typical nonconservative system of two equations \cite{Colombeau4}.
\begin{equation} u_t+(\frac{u^2}{2})_x"="\sigma_x,\end{equation}
\begin{equation} \sigma_t+u\sigma_x"="u_x,\end{equation}
which is nonconservative because of the term $u\sigma_x$ in (56), where $"="$ means that the meaning of this equality is not really clear (it can be = in $\mathcal{G}$ or $\approx$ or a similar concept). As usual we search solutions $u,\sigma$ in form of two constant values $(u_l,\sigma_l)$ on the left and $(u_r,\sigma_r)$ on the right separated by a discontinuity that moves at constant speed $c$, i.e. of the form
\begin{equation} u(x,t)=u_l+[u] H_u(x-ct),\end{equation}
\begin{equation} \sigma(x,t)=\sigma_l+[\sigma] H_\sigma(x-ct),\end{equation}
where $[u]=u_r-u_l, [\sigma]=\sigma_r-\sigma_l$ are respectively the jumps of $u,\sigma$ on the discontinuity located at $x=ct$. Note the basic fact that  we seek a solution in $(\mathcal{G}(\mathbb{R}))^2$ because the term $u\sigma_x$ does not make sense in distribution theory, and that one knows  that one should a priori represent different physical variables by different Heaviside functions, $H_u,H_\sigma$.\\

$\bullet$ \textbf{No shock wave solutions with equality in $\mathcal{G}$}. Using  the equality in $\mathcal{G}(\mathbb{R})$ for both equations (55, 56) we are going to discover they have no shock wave solution of the form (57, 58).\\
\\
 Inserting  (57, 58) into (55) stated with = in $\mathcal{G}(\mathbb{R})$ one obtains
\begin{equation} (-c+u_l)[u]H_u'+[u]^2H_uH_u'=[\sigma]H_\sigma'.\end{equation}
Integration gives
 \begin{equation}H_\sigma=\frac{1}{[\sigma]}((-c+u_l)[u]H_u+\frac{[u]^2}{2}(H_u)^2) + \ constant\end{equation}
 in $\mathcal{G}(\mathbb{R})$. Since $H_u(-\infty)=H_\sigma(-\infty)=0$ it follows that $constant \ =0$. Now $H_u(+\infty)=H_\sigma(+\infty)=1$ gives a formula which is exactly the classical jump formula of the conservative equation (55) (recall that the  jump condition of the equation $u_t+(f(u,v))_x=0$ is $c=\frac{[f]}{[u]}$ where $[f]=f(u_r,v_r)-f(u_l,v_l)$).\\

 Inserting (57, 58) into (56)   stated with = in $\mathcal{G}(\mathbb{R})$   gives  the relation
\begin{equation} (-c+u_l)[\sigma]H_\sigma'+[u][\sigma]H_uH_\sigma'=[u]H_u'
\end{equation}
 between $H_u$ and $H_\sigma$,  very different from (59) since it cannot be integrated like (59) because of the term $H_uH_\sigma'$. One observes that relations (59)  and (61), which follow  from equations (55), (56) respectively, stated with = in $\mathcal{G}$, are incompatible. Indeed these two relations give $H_{\sigma}$ as two different functions of $H_u$. These functions  are respectively (top and bottom of p. 63 in \cite{Colombeau4}):
 \begin{equation} H_\sigma=(1-\frac{[u]^2}{2[\sigma]})H_u+\frac{[u]^2}{2[\sigma]}(H_u)^2, \ \ H_\sigma=\frac{1}{\sigma}Log|1+\frac{2H_u}{-1+2\frac{[\sigma]}{[u]^2}}|.\end{equation}

\textbf{Conclusion: system (55, 56) has no solution of the form (57, 58) when both equations in system (55, 56) are stated with equality in $\mathcal{G}(\mathbb{R})$. This is of course completely unsatisfactory since solutions of the form (57, 58) do exist physically and are even of the greatest importance in physics and engineering for nonconservative systems such as those given above.}\\

$\bullet$\textbf{ Undefinite jump conditions with association}.
Consider system (55, 56) stated in the form of association for both equations and insert (57, 58). From equation (55) one obtains (59) stated in the association sense: it can be integrated and it gives the classical jump formula for (55). From equation (56) one obtains (61) with association in the place of equality, i.e.
\begin{equation} (-c+u_l)[\sigma](H_\sigma)'+[u][\sigma]H_u(H_\sigma)'\approx [u](H_u)'.\end{equation} 
Since this equation is stated in the sense of association it cannot be manipulated by nonlinear calculations  (17): nonlinear calculations  would request the equality in $\mathcal{G}(\mathbb{R})$. Further, as it is, relation (63) cannot be integrated because of the term $H_uH_\sigma'$: $\int H_uH_\sigma'$ can take a continuum of values  depending on the relation between $H_u$ and $H_\sigma$: for instance if $H_u=(H_\sigma)^n$ one obtains $H_uH_\sigma'=(H_\sigma)^nH_\sigma'=\frac{1}{n+1}(H_\sigma^{n+1})'\approx \frac{1}{n+1}H_\sigma'\approx \frac{1}{n+1}\delta$ since $H^n\approx H$ if H is an Heaviside function (equal to 0 to 0 for $x<0$ and to 1 for $x>0$), and  since from   (16) one can differentiate the association.\\

 \textbf{An algebraic relation between $H_u$ and $H_\sigma$, such as for instance one of the relations (62), would be needed to resolve the ambiguity in the product $H_uH_\sigma'$ and to obtain a jump condition from (63)}. \\

We have observed that any such algebraic relation can be obtained by stating (55) or (56) in the sense of equality in $\mathcal{G}(\mathbb{R})$, see (62).\\

Equation (63) shows that $H_u(H_\sigma)'$ is associated to a multiple of the Dirac delta function  ($\delta\approx (H_u)'\approx (H_\sigma)'$). Setting 
\begin{equation} H_u(H_\sigma)'\approx A\delta, \ \ A\in \mathbb{R},\end{equation}
(63) can be integrated and gives
\begin{equation} (-c+u_l)[\sigma]+[u][\sigma]A= [u].\end{equation} 
Now this jump condition depends on the arbitrary parameter $A$ which ca take a continuum of values and contains the ambiguity in the product $H_u(H_\sigma)'$.\\

From natural numerical schemes for equations (55, 56) formulas (64, 65) have been checked with precision: the value $A$ can be computed  from (64) when one has the set of values taken by $H_u$ and $H_\sigma$ in their jumps, which makes sense  since the jumps take place on several cells, see \cite{CLRNP} and \cite{Colombeau4} p. 44-48. One has observed that formula (65) was verified:  \cite{CLRNP} and \cite{Colombeau4} p. 44-48.\\

$\bullet$\textbf{ Definite jump conditions with a mixed strong-weak statement of the equations}.
One has observed on standard systems of  physics \cite{Colombeau4} that when one has a family of $n$ independent physical variables satisfying a system of $n$ equations with no independent equation, the statement of $n-1$ equations with equality in $\mathcal{G}(\mathbb{R}) $ gives at least implicitely all the Heaviside functions representing the jump of $n-1$ variables (by formulas such as one of the formulas (62)) as a function of the Heaviside function of the $n^{th}$ variable, thus resolving the ambiguities in  nonconservative terms. The last Heaviside function remains undetermined among all Heaviside functions in $\mathcal{G}(\mathbb{R})$: this does not matter to obtain the jump conditions: the knowledge of one Heaviside function $H_u$ or $H_\sigma$ in (64) as a function of the other permits the calculation of the number $A$. \\

The jump conditions for the two formulations
\begin{equation}u_t +uu_x=\sigma_x, \ \ \sigma_t+u\sigma_x\approx k^2 u_x, \end{equation}
and
\begin{equation}u_t +uu_x\approx\sigma_x, \ \ \sigma_t+u\sigma_x= k^2 u_x, \end{equation}
$k>0$, have been computed in \cite{Colombeau4} p. 62-64. Jump conditions for systems of $n$ equations, when $n-1$ are stated with $=$ and the last one is stated with $\approx$, have been computed in \cite{Colombeau4} p. 70-77.\\

Numerical investigations to resolve systems stated with = for some equations and $\approx$ for the other equations are reported in \cite{Colombeau4} p. 64-68, \cite{Zalzali} using a large space step for the equations stated with $\approx$ and a small space step for the equations stated with =. One has observed (figure 3) the expected result that the jump conditions are those calculated by stating with = in $\mathcal{G}$ the equation discretized  with small steps and by stating with $\approx$ the equation discretized with large steps. Other numerical investigations consist in the use of a large viscosity for  the equations stated with $\approx$ and a small viscosity for the equations stated with =, (figure 4) \cite{Prignet}. One obtains again the expected result similar to the above for the double scale method.\\

$\bullet$\textbf{Choice of the statement of the equations on physical ground.} To apply the above method for a system of equations modelling a given physical situation one has to decide on physical ground which equations should be stated with = and with $\approx$. This can be delicate since most equations involve approximations. In general one distinguishes between the basic laws (conservation of mass, momentum, total energy  and electric charge) and the state laws (also called constitutive equations) which are extracted from particular experiments and are usually far less accurate. \\

In case of a system of $n$ equations with only one state law, it is natural from physics to state the
state law  with $\approx$, \cite{Colombeau4} p. 69-77. \\

In  the system (42-46) of primordial cosmology there are two state laws: an extension of Newton's law of gravitation (71), and the state law (72) for the energy. It is natural to choose to state (72) with association since it is an approximation.
 \begin {equation} \rho_t+((\rho+\frac{p}{c^2}) u)_x+((\rho+\frac{p}{c^2}) v)_y=0,\end{equation}
\begin {equation}(\rho+\frac{p}{c^2}) ( u_t+uu_x+uv_y)+ p_x+(\rho+\frac{p}{c^2})\Phi_x=0,\end{equation}
\begin {equation}(\rho+\frac{p}{c^2}) ( v_t+vu_x+vv_y)+ p_y+(\rho+\frac{p}{c^2})\Phi_y=0,\end{equation}
\begin {equation}\Delta \Phi=4\pi G(\rho+\frac{3p}{c^2})\end{equation}
\begin{equation}  p\approx\frac{1}{3}\rho c^2 .\end{equation}
This statement would justify the fact cosmologists have already calculated on equations (68-70) to obtain them from the equations of relativity. \\

 For the multifluid system (47-54) in the case of a mixture of a liquid and a gas one naturally states the state law of the liquid with = and the state law of the gas with $\approx$ since the state law of the liquid (which is nearly uncompressible) is far less affected  than the state law of the gas inside the shock waves, \cite {Colombeau4, Berger, AragonaCJ2}.\\

 For the elasticity system (34-41) the situation is more difficult due to the presence of the 3 scalar  state laws (39-41), in fact a 2-D formulation of the Hooke'law which is a tensorial law; an idea is, if the solid is isotropic, to state that all components of the stress are modeled by the same Heaviside function, thus giving three relations between Heaviside functions. Then the statement of the state law (41) with association permits to propose well defined jump conditions.\\

The paths in \cite{DalMaso, LeFloch} are exactly the microscopic profile of Heaviside functions in \cite{CLRNP,Bialivre, Colombeau4,ColombeauJMP}, i.e. the passage between the value 0 for $x<0$ and the value 1 for $x>0$. In this way quantities such as $A$ in (64) are expressed via the paths. The equivalence of the two formulations in the setting of bounded variation functions is proved in \cite{ColombeauHeibig}. The path formulation as it is done in \cite{DalMaso, LeFloch} does not include the concepts of = and $\approx$ which has been  used above to solve the ambiguity in values such as $A$ in (65) due to the products of a Heaviside function and a Dirac distribution, and permits a choice  of  paths in a given physical situation.\\

$\bullet$ \textbf{The continuity equations and the Euler equations in presence of shock waves}. The classical constitutive equation (29) and the classical Euler equation (30), when both stated with equality  in $\mathcal{G}(\mathbb{R})$, imply a remarkable connection between the Heaviside functions $H_{\frac{1}{\rho}}, H_u, H_p$ that depict $\frac{1}{\rho}, u, p$ respectively, according to
\begin{equation} \omega(x,t)=\omega_l+[\omega]H_{\omega}(x-ct),\end {equation} 
$\omega=\frac{1}{\rho}, u, p$.\\

\textbf{Theorem.} \textit{ Equations (29, 30) stated with = in  $\mathcal{G}(\mathbb{R})$ imply equality in  $\mathcal{G}(\mathbb{R})$ of the three a priori different Heaviside functions  $H_{\frac{1}{\rho}},  H_u$ and  $H_p$.}\\

The proof can be found in \cite{Colombeau4} p. 69-72.  Now consider the system of isothermal fluids stated as 
\begin {equation} \rho_t+(\rho u)_x=0,\end{equation}
\begin {equation}(\rho u)_t+(\rho u^2)_x+ p_x=0,\end{equation}
\begin {equation}p\approx C. \rho, \ \ C \ a \ constant.\end{equation}
From above one has  $H_{\frac{1}{\rho}}= H_u= H_p$. We deduce that there are no shock waves that satisfy the state law (76) with equality in $\mathcal{G}(\mathbb{R})$: indeed $p=C\rho$ would imply $H_p=H_\rho$,
therefore  $H_{\frac{1}{\rho}}= H_\rho$. But this last equality is impossible: denoting this joint Heaviside function by $H$ it would imply $(\rho_l+[\rho]H)(\frac{1}{\rho_l}+[\frac{1}{\rho}]H)=1$; using $[\rho]=\rho_r-\rho_l$ and $[\frac{1}{\rho}]=\frac{1}{\rho_r}-\frac{1}{\rho_l}$ one obtains at once $H^2=H$ in $\mathcal{G}(\mathbb{R})$, which is proved to be impossible at the beginning of the paper.\\

$\bullet$ \textbf{The form of the state laws  inside the shock waves.}
This implies that the state law $p=C\rho$ of isothermal gases cannot be satisfied inside the width of the shock waves, although the weak statement (76) implies that it is satisfied on both sides from the formulas $p_l=C \rho_l, p_r=C \rho_r$. Is it possible to modify the statement of (76) so as to state the modified state law with equality in $\mathcal{G}(\mathbb{R})$ inside a shock wave?. Then  \begin {equation}
C=\frac{p}{\rho}=(p_l+[p]H)(\frac{1}{\rho_l}+[\frac{1}{\rho}]H)\end{equation}
is no longer a constant as a generalized function $\in\mathcal{G}(\mathbb{R})$: it varies inside the shock wave. Since $H(x)=\frac{u(x,t)-u_l}{u_r-u_l}=\frac{p(x,t)-p_l}{p_r-p_l}=\frac{\frac{1}{\rho}(x,t)-\frac{1}{\rho}_l}{\frac{1}{\rho}_r-\frac{1}{\rho}_l}$ the function $C=C(x)$ imposed by equality in the state law can be expressed in function of the physical variables $u,p,\rho$  and their left-right values on the shock wave. Then one has at the same time equations (74), (75) and (76) stated with the (strong) equality in $\mathcal{G}$ and presence of shock wave solutions. This  strong statement is suitable to have  uniqueness since strong solutions are  unique in general. \\

$\bullet$ \textbf{A physical interpretation of = and $\approx$ in case of shock waves.} The physical intuitive idea behind the formulation of shock waves with = and $\approx$ is that shock waves would have an infinitesimal width, which is the width of the Heaviside functions used to model the various physical variables (this width has been observed experimentally). This  infinitesimal width could be thought to be of the order of magnitude of a few hundred mean free paths in fluids (a mean free path is the average distance between two collisions) or cristal sizes in solids. This size of the widths of shock waves  is large enough so as to imagine small bags inside these widths, and, by thought, one could state the conservation laws for the volume in these bags (between times $t$ and $t+dt$), which reproduces the classical way physicists obtain the equations of continuum mechanics. This justifies that these  equations should be
\\
\\
\\
\\
\\
\\
\\
\\
\\
\\
\\
\\
\\
\\
\\
\\
\\
\\
\\
\\
\\
\\
\\
\\
\\
\\
\\

\begin{figure}[h]
\centering
\includepdf[width=\textwidth]{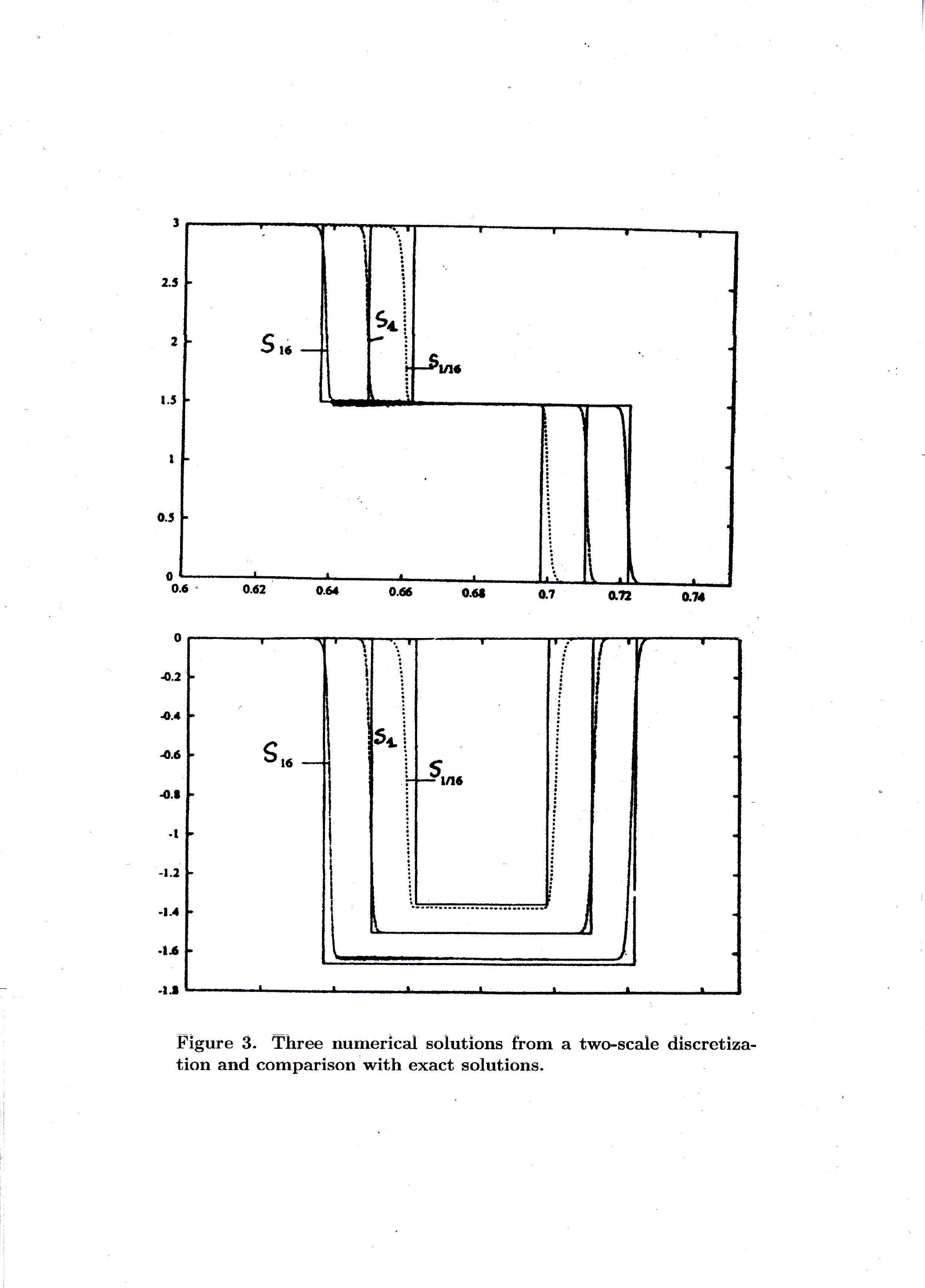}
\end{figure}
\textit{}
\\
stated with the (strong) equality as done above. On the other hand the state laws stated with $\approx$ are only valid on both sides of the shock waves, in the same way as $H$ and $H^2$ are equal on both sides of the discontinuity, and are different inside the discontinuity. By definition of a state law we know that such a law  has been checked (up to now) by experimental physicists in the case of (quasi) absence of variation of the physical variables (absence of experimental apparatus to observe inside widths of shock waves), that is, on both sides of the shock waves but not  inside the width of the shock waves. An interesting fact is that the genuine state law that holds inside the width of the shock waves can be calculated as a consequence of the equations when all are stated with the strong equality (77), knowing the state laws outside of the shock waves.\\

\vskip 25cm
\begin{figure}[h]
\centering
\includepdf[width=\textwidth]{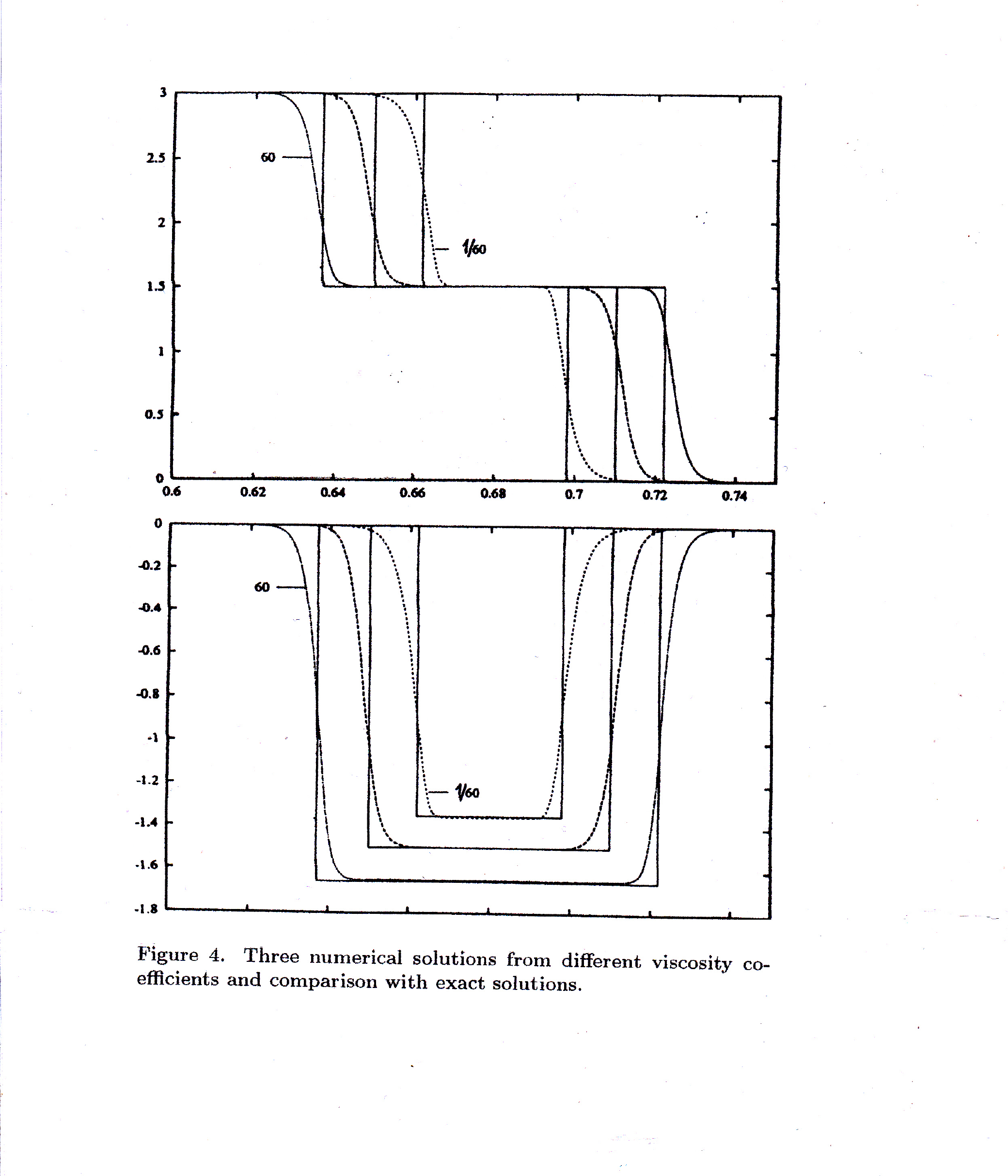}
\end{figure}
$\bullet$ \textbf{Numerical investigations.}
Numerical tests on the statement of equations with = and $\approx$ are reported below (figures 3 and 4).  Various numerical methods have been devised to treat systems in which some equations are stated with $=$ and the other ones with $\approx$. These methods consist in a discretization of the equations stated with $\approx$ which should be much coarser than the discretization of the equations stated with =. In \cite{Zalzali} one uses two discretization steps: a large step for the equations stated with $\approx$ and a small one for the equations stated with =. One obtains the expected results: one can observe that the numerical solutions so obtained coincide with the explicit jump conditions corresponding to the various statements (figure 3). In \cite{Prignet} one uses viscous second members with  viscosity coefficients larger for the equations stated with $\approx$. Again one obtains the expected results (figure 4). \\

 In figures 3 and 4 we present the exact solutions for the system (55, 56) stated successively with $(\approx,=), (\approx, \approx)$ with same Heaviside functions $H_u$ and $H_\sigma$, and $(=,\approx)$, and numerical results from various techniques of discretization of the equations.\\

 In figure 3 one uses a two scale method for the discretization of derivatives: for equation $ u_t+uu_x "="\sigma_x$ one sets $f'=\frac{f(x,t)-f(x-h_1,t)}{h_1}, \ f=u,\sigma$ and for equation $\sigma_t+u\sigma_x"="u_x$ one sets $f'=\frac{f(x,t)-f(x-h_2,t)}{h_2}, \ f=u,\sigma$. One tests various values of $k:=\frac{h_1}{h_2}$. For $k=16$ the first equation is discretized in a far coarser way than the second and we observe the numerical result coincides with the statement of the first equation with $\approx$ and the second with = , as it could be expected. For $k=1$ one observes the numerical result corresponds to the same Heaviside functions for $u$ and $\sigma$. For $k=\frac{1}{16}$ we observe the numerical result corresponds to the statement of the first equation with = and the second with $\approx$ as expected.\\

 In figure 4 the system is stated with a small viscosity i.e.

$$u_t+uu_x=\sigma_x+\epsilon_1u_{xx},$$   $$\sigma_t+u\sigma_x=u_x+\epsilon_2\sigma_{xx},$$ 
\\
and one tests the quotient $k:=\frac{\epsilon_2}{\epsilon_1}$. As expected, for $k=\frac{1}{60}$ the fact that $\epsilon_2<<\epsilon_1$ implies that one observes the solution $(\approx,=)$. For $k=1$ one observes same Heaviside functions and for $k=60$ one observes the solution $ (=, \approx)$. This shows that the various statements of equations with = and $\approx$ are a fact of the real world (here numerical calculations) and that it can be mastered numerically.\\
In \cite{Joseph} one proves that the above viscous Cauchy problem admits solutions in the association sense when $k=1$ for rather arbitrary initial data in $\mathcal{G}(\mathbb{R})$. This is done by solving the equations for given $\epsilon>0$. The result is unknown when $k\not=1$.\\

\textbf{4. A topological continuation of the nonlinear generalized functions: S\~ao Paulo 2012.} In the theory of Banach spaces compactness imposes the use of weak topologies but these topologies are uncompatible with nonlinearity in that $u_n\rightarrow u, v_n\rightarrow v$ do not imply in general $u_nv_n\rightarrow uv.$ Correspondingly  the locally convex spaces of distributions  are not topological algebras. The aim of this section is to present a context containing irregular objects such as all distributions with compact support in which one has compactness and nonlinearity at the same time.  We construct subalgebras of $\mathcal{G}(\Omega)$ that enjoy Hausdorff locally convex topologies  suitable for the development of a  functional analysis: they have very good topological properties in particular concerning compactness (strong duals of Fr\'echet-Schwartz spaces) at the same time as they are compatible both with partial derivatives and nonlinearity.\\

In  \cite{AragonaCJ} we construct a family of differential algebras generically denoted $\mathcal{G}_a(\Omega)$  in the situation $\mathcal{G}_a(\Omega)\subset \mathcal{G}_s(\Omega)$ such that\\
\\
$\bullet$ The multiplication $ \mathcal{G}_a(\Omega)\times\mathcal{G}_a(\Omega)\longmapsto\mathcal{G}_a(\Omega), (F,G)\longmapsto FG$, and the differentiation $\mathcal{G}_a(\Omega)\longmapsto\mathcal{G}_a(\Omega), F\longmapsto \frac{\partial^{|\alpha|} F}{\partial x_1^{\alpha_1},\dots,\partial x_n^{\alpha_n}}$, are continuous,\\
\\
$\bullet$  $\mathcal{G}_a(\Omega)$ is a Hausdorff locally convex topological algebra which is a strong dual of a Frechet-Schwartz space (in particular for any bounded set $b$ there is a larger bounded set $B$ which is convex and balanced such that the vector space $\cup_n nB$ normed with the gauge of $B$ is a Banach space and such that  the set $b$ is relatively compact in this  Banach space),\\
\\
 $\bullet$  spaces of distributions such as $\mathcal{E}'(\Omega)$ are contained in  $\mathcal{G}_a(\Omega)$ with continuous inclusions so that any bounded set in these spaces of distributions is bounded in $\mathcal{G}_a(\Omega)$, therefore relatively compact as explained above.\\

The construction in \cite{AragonaCJ} is based on the properties of holomorphic functions: uniqueness of analytic continuation, classification of point singularities, Laurent expansion around a pole. This permits to have only one representative in each class of the quotient $\frac{reservoir}{kernel}$ defining $\mathcal{G}_s(\Omega)$. The absence of a quotient so obtained permits to have nice locally convex topologies. \\

This could be an interesting continuation of the  theory of nonlinear generalized functions started inside the Leopoldo Nachbin group 30 years ago, in which a topological structure with such properties remained lacking up to now. A topology introduced early  by H.A. Biagioni \cite{Bialivre} (but somewhat hidden in a minor point p. 44-45) and rediscovered  independently by D. Scarpalezos \cite{Scarpatopo, Scarpatopo2} is not a locally convex vector space topology and has no known compactness property. Anyway it has interested applications: \cite{AragonaFernandez, AragonaFerJu1, Novo,  AragonaFerJu2,  AragonaFerJu3, AragonaGarcia, Bialivre, Garetto, GarettoV, Vernaeve,Scarpatopo, Scarpatopo2} .\\

\end{document}